# On the computation of the $n^{th}$ decimal digit of various transcendental numbers.

by Simon Plouffe

November 30, 1996
Revised December 1, 2009

## Abstract

A method for computing the $n^{th}$ decimal digit of $\pi$ in $O(n^3 log(n)^3)$ time and with very little memory is presented here. The computation is based on the recently discovered Bailey-Borwein-Plouffe algorithm and the use of a new algorithm that simply splits an ordinary fraction into its components. The algorithm can be used to compute other numbers like (3), $\pi\sqrt{3}$, $\pi^2$ and $\frac{2}{\sqrt{5}} ln(\varphi)$ where $\varphi$ is the golden ratio. The computation can be achieved without having to compute the preceding digits. I claim that the algorithm has a more theoretical rather than practical interest, I have not found a faster algorithm or proved that it cannot exist.

The formula for $\pi$ used is

$$\pi + 3 = \sum_{n=1}^{\infty} \frac{n 2^n n!^2}{(2n)!}$$

Introduction
Key observation and the Splitting Algorithm
Other numbers
Conclusion and later developments
Bibliography

## Introduction

The computation of the $n^{th}$ digit of irrational or transcendental numbers was considered either impossible or as difficult to compute the number itself. In 1995 [BBP], have found a new way of computing the $n^{th}$ binary digit of various constants like $\pi$ and $ln(2)$.

An intensive computer search was then carried out to find if that algorithm could be used to compute a number in an arbitrary base. I present here a way of computing the $n^{th}$ decimal digit of $\pi$ (or any other base) by using more time than the [BBP] algorithm but still with very little memory.

## Key observation and formula

The observation is that a fraction $1/ab$ can be split into $k_1/a + k_2/b$ by using the continued fraction algorithm of a/b. Here a and b are two prime powers. This is equivalent to having to solve a Diophantine equation for $k_1$ and $k_2$ - it is always possible to do so if (a,b) = 1, see [HW] if they have no common factor. If we have more than 2 prime factors then it can be done by doing 2 at the time and then using the result to combine with the third element. This way an arbitrary big integer M can be split into small elements. If we impose the conditions on M of having only small factors (meaning that the biggest prime power size is of the order of a computer word), then an arbitrary M can be represented. If this is true then a number of known series and numbers can then be evaluated. For example the expression $1/\binom{2n}{n}$ the central binomials satisfy that: the prime powers of this number are small when n is big.

Example:

$$1/\binom{100}{50} = \frac{1}{100891344545564193334812497256} =$$

$$\frac{1}{2^3 \cdot 3^4 \cdot 11 \cdot 13 \cdot 17 \cdot 19 \cdot 29 \cdot 31 \cdot 53 \cdot 59 \cdot 61 \cdot 67 \cdot 71 \cdot 73 \cdot 79 \cdot 83 \cdot 89 \cdot 97}$$

Now if we take 2 elements at the time and solve the simple Diophantine equation and proceed this way:

1) $\frac{1}{ab} = \frac{k_1}{a} + \frac{k_2}{b}$

2) $\left(\frac{\frac{k_1}{a}+\frac{k_2}{b}}{c}\right) = \frac{m_1}{a} + \frac{m_2}{b} + \frac{m_3}{c}$

3) We proceed with the next element.

At each step the constants $k_1$ and $k_2$ are determined by simply expanding a/b into a continued fraction and keeping the *before last* continuant, later $m_1$, $m_2$ and $m_3$ are determined the same way. Having finished with that number we quickly arrive at a number which is (modulo 1) the same number but represented as a sum of only small fractions.

$$\frac{1}{100891344545564193334812497256} =$$

$$\frac{5}{8} + \frac{20}{81} + \frac{10}{11} + \frac{2}{13} + \frac{13}{17} + \frac{10}{19} + \frac{4}{29} + \frac{5}{31} + \frac{23}{53} + \frac{41}{59} + \frac{29}{61} + \frac{37}{67} + \frac{33}{71} + \frac{19}{73} + \frac{36}{79} + \frac{7}{83} + \frac{13}{89} + \frac{88}{97}$$

The time taken to compute this expression is $O(log(n)n^2)$, $log$(n) being the time spent to compute with the Euclidian algorithm on each number. I did not take into account the time spent on finding what is the next prime in the expression simply because we can consider (at least for the moment) that the applicability of the algorithm is a few thousands digits and so the time to compute a prime is really a matter of a few seconds in that range for the whole process. Since we know by advance what is the maximal prime there could be in $\binom{2n}{n}$ then we can do it with a greedy algorithm that pulls out the factors until we reach 2n, and this can be done without having to compute the actual number which would obviously not fit into a small space. It can be part of the loop without having to store any number apart from the current n. For any p in $\binom{2n}{n}$ the maximal exponent is (as Robert Israel pointed out).

$$\sum_{k=1}^{2n}\left(\left\lceil\frac{2n}{p^k}\right\rceil - \left\lfloor\frac{n}{p^k}\right\rfloor\right)$$

Equivalently, for p = 2 it gives the number of '1' in the binary expansion of n, for p = 3 there is another clue with the ternary expansion of the number and the number of times the pattern '12' appears. Now looking at the $\sum_{n=1}^{\infty}\frac{1}{\binom{2n}{n}}$ term we can say that the series is essentially $\pi\sqrt{3}$ since it differs only by $\frac{4}{9}\pi\sqrt{3}+\frac{1}{3}$ since these are 2 small rational numbers we can use BBP algorithm to carry the computation to an arbitrary position in almost no time. Having $n/\binom{2n}{n}$ instead of (1) only simplifies the process.

To compute the final result of each term we need only few memory elements,

1 for the partial sum so far. (evaluated later with the BBP algorithm). 4 for the current fractions $k_1$/a and $k_2$/b. 2 for the next element to be evaluated: 1/c. 1 for n itself.

So with as little as 8 memory elements the sum for each term of (1) can be carried out a without having to store any number greater than a computer word in log(n) time, adding this for each element the total cost for (1) is then $n^3 log(n)$.

The next thing we have to consider is that , if we have an arbitrary large M and if M has only small factors then 1/M can be computed. First, we need to represent 1/M as $\sum_{j=1}^{k}\frac{a_i}{p_i^j}$ where $p_i^j$ is a prime power and $a_i$ is smaller than $p_i^j$.

If we have $2^n$/M then by using the binary method on each element of the representation of 1/M with (2) is possible in $log(n)$ time. Again if we don't want to store the element of (2) in memory we can do it as we do the computation of
the first part at each step. In this algorithm we can either store the powers of 2 to do the binary method or not. There is variety of ways to do it, we refer to [Knuth vol. 2] for explanations.

This step is important, essentially once we can represent 1/(a*b) by splitting them then to multiply by $2^n$ only adds $log(n)$ steps for each element and it can be done in arbitrary base since we have the actual fraction for each element of (2). It only pushes the decimal (or the decimal point of the base chosen), further. At any moment only one element in the expansion of 1/M is considered with the current fraction, that same fraction can be represented in base 10 at any time if we want the decimal expansion at that point. For this reason multiplying the current fraction by $2^n$ involves only small numbers and fractions.

Once this is done, the total cost becomes $n^3\log(n)^2$. This cost is for the computation of the $k^{th}$ partial sum of

$$\sum_{n=1}^{\infty} \frac{2^n}{n \binom{2n}{n}}$$

where $\binom{2n}{n}$ is the central binomial coefficient. If we want at each step to compute (the final $n^{th}$ digit) then we need $log(n)$ steps to do it. It can be done in any base chosen in advance, in BBP the computation could be done in base 2 but here we have the actual explicit fraction which is independent of the base. This is where we actually compute the decimal expansion of the final fraction of the process. So finally the $n^{th}$ digit of $\pi$ can be computed in $n^3\log(n)^3$ steps.

Other Numbers

By looking at the plethora of formulas of the same type as (1) or (3) we see that [RamI and IV] the numbers $\pi\sqrt{3}$, $\pi$, $\zeta(3)$ and even powers of $\pi$ can be computed as well. The condition we need to ensure is: if any term of a series can be split into small fractions of size no greater than that of a computer word, then it is part of that class. This includes series of the type:

$$\sum_{n=1}^{\infty} \frac{c^n}{P(n) \binom{2mn}{n}^r}$$

where c is an integer, P(n) a polynomial and $\binom{2mn}{n}$ is a near central binomial coefficient. This class of series contains many numbers that are not yet identified in terms of known constants and conversely the known constants that are of similar nature like $\zeta(5)$ have not yet been identified as members of the class. The process of identifying a series as being expressed in terms of known constants and the exact reverse process is what the Inverter tries to do.

The number e or exp(1) which is $\sum_{n=1}^{\infty} \frac{1}{n!}$ does not satisfies our condition because 1/n! eventually contains high powers of 2, therefore cannot be computed to the $n^{th}$ digit using our algorithm. The factorisation of 1/n! has high powers of small primes, the highest is $2^k$ and k is nearly the size of n. For this particular number only very few series are known and appear to be only a variation on that first one.

Others like gamma or Catalan do not seem to have a proper series representation and computer search using Ferguson's PSLQ or LLL with Maple and Pari-Gp gave no answer to this. Algebraic numbers like $\sqrt{2}$ have not been yet been fully investigated and we still do not know if those would fall into this category.

## Conclusions and later developments

There are many, but first and foremost we cannot resist thinking at William Shanks who did the computation of $\pi$ by hand in 1853 - if he would had known this algorithm, he would have certainly tried it before spending 20 years of his life computing $\pi$ (half of it on a mistake). Secondly, the algorithm shown here is theoretical and not practical. We do not know if there is a way to improve it, and if so then it is reasonable to think that it could then be used to check long computations like the one that Daisuke Takahashi conducted in August 2009 for the computation of $\pi$ to 2576 billion digits. There could be a way to speed the algorithm to make it an efficient algorithm.

Thirdly, so far there are 2 classes of numbers that can be computed to the $n^{th}$ digit:

The SC(2) class as in the [BBP] algorithm which includes various polylogarithms.

This new class of numbers. Now what's next ?, so far we do not know whether, for example, series whose general term is H(n)/$2^n$ (where H(n) is the $n^{th}$ harmonic number) which fall into the first class, can be extended. We think that this new approach is only the tip of the iceberg. Finally, it is interesting to observe that we can then compute $\pi$ to the 1000000'th digit without having to store (hardly) any array or matrix, so it can be computed using a small pocket calculator. We also note that, in some way we have a way to produce the digits of Pi without using memory, this means that the number is compressible , if we consider that we could use the algorithm to produce a few thousands digits of the number. We think that other numbers are yet to come and that there is a possibility (?) of having a direct formula for the n'th digit (in any base) of a naturally occurring constant like $ln(2)$.

Fabrice Bellard improved this algorithm in 1997 to $O(n^2)$ as explained here: http://bellard.org/pi/pi_n2/pi_n2.html

Xavier Gourdon made also an improvement here: http://numbers.computation.free.fr/Constants/Algorithms/nthdecimaldigit.pdf

This later algorithm is not based on the same idea but inspired from it. It uses a series of $\pi$ plus the asymptotic representation of it. When combined it results an algorithm that can reach 4000000 digits in a matter of hours.

## Acknowledgments

I wish to thanks, Robert Israel (Univ. British Columbia) and David H. Bailey [NASA] for their helpful comments.

# Bibliography

[BBP] David H. Bailey, Péter B. Borwein and Simon Plouffe, On The Rapid Computation of Various Polylogarithmic Constants, april 1997 in Mathematics of Computation.

[RAM] Bruce C. Berndt, Ramanujan Notebooks vols. I to V. Springer Verlag, New York.

[RI] Robert Israel at University of British Columbia, personal communication.

[AS] M. Abramowitz and I. Stegun, Handbook of Mathematical Functions, Dover, New York, 1964.

[HW] G. H. Hardy and E. M. Wright, An Introduction to the Theory of Numbers 5e, Oxford University Press, 1979.

[Shanks] W. Shanks, Contributions to Mathematics Comprising Chiefly of the Rectification of the Circle to 607 Places of Decimals, G. Bell, London, 1853.

[Knuth] D.E. Knuth, The Art of Computer Programming, Vol. 2: Seminumerical Algorithms, Addison-Wesley, Reading, MA, 1981.

[PI] Plouffe Inverter at http://pi.lacim.uqam.ca/eng Inverseur de Plouffe.

[Bellard] F. Bellard, Computation of the n'th digit of pi in any base in O(n2), unpublished (1997) http://fabrice.bellard.free.fr/pi/pi n2/pi n2.html

[Gourdon] Xavier Gourdon, Computation of the n-th decimal digit of _ with low memory available here: February 11, 2003.
http://numbers.computation.free.fr/Constants/Algorithms/nthdecimaldigit.pdf